\documentclass{amsart}
\usepackage{graphicx,amscd, color,amsmath,amsfonts,amssymb,geometry,amssymb}
\providecommand{\U}[1]{\protect\rule{.1in}{.1in}}
\newtheorem{theorem}{Theorem}[section]

\newtheorem{remark}[theorem]{Remark}

\newtheorem{lemma}[theorem]{Lemma}

\numberwithin{equation}{section}

\begin{document}
\title{New upper bounds for the constants in the Bohnenblust-Hille inequality}
\author[D. Pellegrino \and J. B. Seoane-Sep\'{u}lveda]{Daniel Pellegrino{*} \and Juan B. Seoane-Sep\'{u}lveda\textsuperscript{**}}

\address{Departamento de An\'{a}lisis Matem\'{a}tico,\newline\indent Facultad de Ciencias Matem\'{a}ticas, \newline\indent Plaza de Ciencias 3, \newline\indent Universidad Complutense de Madrid,\newline\indent Madrid, 28040, Spain.}
\email{jseoane@mat.ucm.es}

\address{Departamento de Matem\'{a}tica, \newline\indent Universidade Federal da Para\'{\i}ba, \newline\indent 58.051-900 - Jo\~{a}o Pessoa, Brazil.} \email{pellegrino@pq.cnpq.br}
\thanks{2010 Mathematicss Subject Classification: 46G25, 47L22, 47H60.}
\thanks{\textsuperscript{**}Supported by CNPq Grant 620108/2008-8 (Edital Casadinho).}
\thanks{\textsuperscript{**}Supported by the Spanish Ministry of Science and
Innovation, grant MTM2009-07848.}
\begin{abstract}
A classical inequality due to Bohnenblust and Hille states that for every
positive integer $m$ there is a constant $C_{m}>0$ so that
$$\left(  \sum\limits_{i_{1},\ldots ,i_{m}=1}^{N}\left\vert U(e_{i_{^{1}}},\ldots ,e_{i_{m}})\right\vert ^{\frac{2m}{m+1}}\right)  ^{\frac{m+1}{2m}}\leq C_{m}\left\Vert U\right\Vert$$
for every positive integer $N$ and every $m$-linear mapping $U:\ell_{\infty}^{N}\times\cdots\times\ell_{\infty}^{N}\rightarrow\mathbb{C}$, where
$C_{m}=m^{\frac{m+1}{2m}}2^{\frac{m-1}{2}}.$ The value of $C_{m}$ was improved
to $C_{m}=2^{\frac{m-1}{2}}$ by S. Kaijser and more recently H. Qu\'{e}ffelec
and A. Defant and P. Sevilla-Peris remarked that $C_{m}=\left(  \frac{2}{\sqrt{\pi}}\right)^{m-1}$ also works. The Bohnenblust--Hille inequality also holds for real Banach spaces with the constants $C_{m}=2^{\frac{m-1}{2}}$.
In this note we show that a recent new proof of the Bohnenblust--Hille
inequality (due to Defant, Popa and Schwarting) provides, in fact, quite better
estimates for $C_{m}$ for all values of $m \in \mathbb{N}$. In particular, we will also show that, for real scalars, if $m$ is even with $2\leq m\leq 24$, then
$$C_{\mathbb{R},m}=2^{\frac{1}{2}}C_{\mathbb{R},m/2}.$$
We will mainly work on a paper by Defant, Popa and Schwarting, giving some remarks about their work and explaining how to, numerically, improve the previously mentioned constants.
\end{abstract}
\maketitle

\section{Preliminaries and background}

In 1930, Littlewood proved that
\[
\left(  \sum\limits_{i,j=1}^{N}\left\vert U(e_{i},e_{j})\right\vert ^{\frac
{4}{3}}\right)  ^{\frac{3}{4}}\leq\sqrt{2}\left\Vert U\right\Vert
\]
for every bilinear form $U:\ell_{\infty}^{N}\times\ell_{\infty}^{N}
\rightarrow\mathbb{C}$ and every positive integer $N.$ This is the well-known
Littlewood's $4/3$ inequality \cite{Litt}.

One year later, in 1931, Bohnenblust and Hille (\cite{bh}) improved this
result to multilinear forms (see also \cite{defant2, defant} for recent
approaches). More precisely, the Bohnenblust--Hille inequality asserts that
for every positive integer $m$ there is a $C_{m}>0$ so that
\[
\left(  \sum\limits_{i_{1},\ldots,i_{m}=1}^{N}\left\vert U(e_{i_{^{1}}},
\ldots,e_{i_{m}})\right\vert ^{\frac{2m}{m+1}}\right) ^{\frac{m+1}{2m}}\leq
C_{m}\left\Vert U\right\Vert
\]
for every $m$-linear mapping $U:\ell_{\infty}^{N}\times\cdots\times
\ell_{\infty}^{N}\rightarrow\mathbb{C}$ and every positive integer $N$ (for
polynomial versions of Bohnenblust--Hille inequality we refer to
\cite{annals2011}). The original upper estimate for $C_{m}$ is $m^{\frac
{m+1}{2m}}2^{\frac{m-1}{2}}$, but several improvements have been obtained
since then. For instance, as an illustration for the complex case we compare,
below, and for some values of $m$, the original constants with the
improvements obtained by S. Kaijser \cite{Ka} and H. Qu\'{e}ffelec \cite{Que},
Defant, P. Sevilla-Peris \cite{defant2}:%

\begin{tabular}
[c]{c|c|c|c}%
$m$ & $C_{m}=\left(  \frac{2}{\sqrt{\pi}}\right)  ^{m-1}${\small (\cite{Que,
defant2}, 1995)} & $C_{m}=2^{\frac{m-1}{2}}${\small (\cite{Ka}, 1978)} &
$C_{m}=m^{\frac{m+1}{2m}}2^{\frac{m-1}{2}}${\small (\cite{bh},1931)}\\\hline
$3$ & $\approx1.273$ & $2$ & $\approx4.160$\\
$4$ & $\approx1.437$ & $\approx2.828$ & $\approx6.726$\\
$5$ & $\approx1.621$ & $4$ & $\approx10.506$\\
$6$ & $\approx1.829$ & $\approx5.657$ & $\approx16.088$\\
$7$ & $\approx2.064$ & $8$ & $\approx24.322$\\
$8$ & $\approx2.330$ & $\approx11.314$ & $\approx36.442$\\
$9$ & $\approx2.628$ & $16$ & $\approx54.232$\\
$10$ & $\approx2.965$ & $\approx22.627$ & $\approx80.283$\\
$15$ & $\approx5.425$ & $128$ & $\approx542.574$\\
$50$ & $\approx372$ & $\approx23,726,566$ & $\approx174,465,512$\\
$100$ & $\approx155,973$ & $\approx7.96131459\cdot10^{14}$ & $\approx
8.14675743\cdot10^{15}$%
\end{tabular}
\bigskip

The Bohnenblust--Hille inequality also holds for real Banach spaces but
sharper estimates for $C_{m}$, in this case, seem to be $C_{m}=2^{\frac
{m-1}{2}}$.

The aim of this paper is show that improved values for $C_{m}$ are essentially
contained in \cite{defant} for both real and complex cases.

The (complex and real) Bohnenblust--Hille inequality can be re-written in the
context of multiple summing multilinear operators, as we will see next.
Multiple summing multilinear mappings between Banach spaces is a recent, very
important and useful nonlinear generalization of the concept of absolutely
summing linear operators. This class was introduced, independently, by Matos
\cite{collect} (under the terminology fully summing multilinear mappings) and
Bombal, P\'{e}rez-Garc\'{\i}a and Villanueva \cite{bpgv}. The interested
reader can also refer to \cite{PAMS2008,REMC2010} for other Bohnenblust--Hille
type results.

Throughout this paper $X_{1},\ldots,X_{m}$ and $Y$ will stand for Banach
spaces over $\mathbb{K}=\mathbb{R}$ or $\mathbb{C}$, and $X^{\prime}$ stands
for the dual of $X$. By $\mathcal{L}(X_{1},\ldots,X_{m};Y)$ we denote the
Banach space of all continuous $m$-linear mappings from $X_{1}\times
\cdots\times X_{m}$ to $Y$ with the usual sup norm. For $x_{1}, \cdots,x_{n} $
in $X$, let%

\[
\Vert(x_{j})_{j=1}^{n}\Vert_{w,1}:=\sup\{\Vert(\varphi(x_{j}))_{j=1}^{n}
\Vert_{1}:\varphi\in X^{\prime},\Vert\varphi\Vert\leq1\}.
\]
If $1\leq p<\infty$, an $m$-linear mapping $U\in\mathcal{L}(X_{1},\ldots
,X_{m};Y)$ is multiple $(p;1)$-summing (denoted $\Pi_{(p;1)}(X_{1}
,\ldots,X_{m};Y)$) if there exists a constant $L_{m}\geq0$ such that
\begin{equation}
\left(  \sum_{j_{1},\ldots,j_{m}=1}^{N}\left\Vert U(x_{j_{1}}^{(1)}
,\ldots,x_{j_{m}}^{(m)})\right\Vert ^{p}\right)  ^{\frac{1}{p}}\leq L_{m}
\prod_{k=1}^{m}\left\Vert (x_{j}^{(k)})_{j=1}^{N}\right\Vert _{w,1}\label{lhs}%
\end{equation}
for every $N\in\mathbb{N}$ and any $x_{j_{k}}^{(k)}\in X_{k}$, $j_{k}
=1,\ldots,N$, $k=1,\ldots,m$. The infimum of the constants satisfying
(\ref{lhs}) is denoted by $\left\Vert U\right\Vert _{\pi(p;1)}$. For $m=1$ we
have the classical concept of absolutely $(p;1)$-summing operators (see, e.g.
\cite{defant3, Di}).

A simple reformulation of the Bohnenblust--Hille inequality asserts that every
continuous $m$-linear form $T:X_{1}\times\cdots\times X_{m}\rightarrow
\mathbb{K}$ is multiple $(\frac{2m}{m+1};1)$-summing with $L_{m}
=C_{m}=2^{\frac{m-1}{2}}$ (or $L_{m}=\left(  \frac{2}{\sqrt{\pi}}\right)
^{m-1}$ for the complex case, using the estimates of Defant and Sevilla-Peris,
\cite{defant2}), although in the real case the best known constants seem to be
$C_{m}=2^{\frac{m-1}{2}}$.

The main goal of this paper is to calculate \emph{better} constants for the
Bohnenblust--Hille inequality in the real and complex case (which are derived
from \cite{defant}). For this task we will explore the proof of a general
vector-valued version of Bohnenblust--Hille inequality (\cite[Theorem
5.1]{defant}).

Let us recall some results that we will need in this note. The first one is a
well-known inequality due to Khinchine (see \cite{Di}):

\begin{theorem}
[Khinchine's inequality]\label{k} For all $0<p<\infty$, there exist constants
$A_{p}$ and $B_{p}$ such that
\begin{equation}
A_{p}\left(  \sum_{n=1}^{N}\left\vert a_{n}\right\vert ^{2}\right)  ^{\frac
{1}{2}}\leq\left(  \int_{0}^{1}\left\vert \sum_{n=1}^{N}a_{n}r_{n}\left(
t\right)  \right\vert ^{p}dt\right)  ^{\frac{1}{p}}\leq B_{p}\left(
\sum_{n=1}^{N}\left\vert a_{n}\right\vert ^{2}\right)  ^{\frac{1}{2}%
}\label{lpo}%
\end{equation}
for every positive integer $N$ and scalars $a_{1}, \ldots,a_{n}$ (here,
$r_{n}$ denotes the $n-th$ Rademacher function)$.$
\end{theorem}

Above, it is clear that $B_{2}=1.$ From (\ref{lpo}) it follows that
\begin{equation}
\left(  \int_{0}^{1}\left\vert \sum_{n=1}^{N}a_{n}r_{n}\left(  t\right)
\right\vert ^{p}dt\right)  ^{\frac{1}{p}}\leq B_{p}A_{r}^{-1}\left(  \int
_{0}^{1}\left\vert \sum_{n=1}^{N}a_{n}r_{n}\left(  t\right)  \right\vert
^{r}dt\right)  ^{\frac{1}{r}}\label{kcc}%
\end{equation}
and the product of the constants $B_{p}A_{r}^{-1}$ will appear later on in
Theorem \ref{d}. The notation $A_{p}$ and $B_{p}$ will be kept along the
paper. Next, let us recall a variation of an inequality due to Blei (see
\cite[Lemma 3.1]{defant}).

\begin{theorem}
[Blei, Defant et al.]\label{b} Let $A$ and $B$ be two finite non-void index
sets, and $(a_{ij})_{(i,j)\in A\times B}$ a scalar matrix with positive
entries, and denote its columns by $\alpha_{j}=(a_{ij})_{i\in A}$ and its rows
by $\beta_{i}=(a_{ij})_{j\in B}.$ Then, for $q,s_{1},s_{2}\geq1$ with
$q>\max(s_{1},s_{2})$ we have
\[
\left(  \sum_{(i,j)\in A\times B}a_{ij}^{w(s_{1},s_{2})}\right)  ^{\frac
{1}{w(s_{1},s_{2})}}\leq\left(  \sum_{i\in A}\left\Vert \beta_{i}\right\Vert
_{q}^{s_{1}}\right)  ^{\frac{f(s_{1},s_{2})}{s_{1}}}\left(  \sum_{j\in
B}\left\Vert \alpha_{j}\right\Vert _{q}^{s_{2}}\right)  ^{\frac{f(s_{2}
,s_{1})}{s_{2}}},
\]
with
\begin{align*}
w  &  :[1,q)^{2}\rightarrow\lbrack0,\infty),\text{ }w(x,y):=\frac
{q^{2}(x+y)-2qxy}{q^{2}-xy},\\
f  &  :[1,q)^{2}\rightarrow\lbrack0,\infty),\text{ }f(x,y):=\frac{q^{2}
x-qxy}{q^{2}(x+y)-2qxy}.
\end{align*}

\end{theorem}

The following theorem is a particular case of \cite[Lemma 2.2]{defant} for
$Y=\mathbb{K}$ using that the cotype $2$ constant of $\mathbb{K}$ is $1$,
i.e., $C_{2}(\mathbb{K})=1$ (following the notation from \cite{defant})$:$

\begin{theorem}
[Defant et al]\label{d} Let $1\leq r\leq2$, and let $(y_{i_{1},\ldots,i_{m}
})_{i_{1},\ldots,i_{m}=1}^{N}$ be a matrix in $\mathbb{K}$. Then
\begin{align*}
&  \left(  \sum\limits_{i_{1},\ldots,i_{m}=1}^{N}\left\vert y_{i_{1}, \ldots,
i_{m} }\right\vert ^{2}\right)  ^{1/2} \leq\left(  A_{2,r}\right) ^{m}\left(
\int\nolimits_{[0,1]^{m}}\left\vert \sum\limits_{i_{1}, \ldots,i_{m}=1}
^{N}r_{i_{1}}(t_{1}) \ldots r_{i_{m}}(t_{m})y_{i_{1} \ldots i_{m}}\right\vert
^{r} dt_{1} \ldots dt_{m}\right) ^{1/r},
\end{align*}
where
\[
A_{2,r}\leq A_{r}^{-1}B_{2}=A_{r}^{-1}\text{ (since }B_{2}=1\text{).}%
\]

\end{theorem}

The meaning of $A_{2,r}$, $w$ and $f$ from the above theorems will also be
kept in the next section and $K_{G}$ will denote the complex Grothendieck
constant. Also, and throughout the paper, $C_{\mathbb{K},m}$ will denote our
estimates on the constants for the real or complex case ($\mathbb{K}%
=\mathbb{R}$ or $\mathbb{C}$ respectively).

\section{Improved constants for the Bohnenblust--Hille theorem: The real case}

From the proof of \cite[Corollary 5.2]{defant} and using the optimal values
for the constants of Khinchine's inequality (due to Haagerup), the following
estimates can be calculated for $C_{m}$:
\begin{align}
C_{\mathbb{R},2}  &  =\sqrt{2},\text{ }\\
C_{\mathbb{R},m}  &  =2^{\frac{m-1}{2m}}\left(  \frac{C_{\mathbb{R},m-1}
}{A_{\frac{2m-2}{m}}}\right)  ^{1-\frac{1}{m}}\text{ for }m\geq3.\label{qqssw}%
\end{align}
In particular, if $2\leq m\leq13$,
\begin{equation}
\text{ }C_{\mathbb{R},m}=2^{\frac{m^{2}+m-2}{4m}}.\label{qq8}%
\end{equation}
If we change a little bit the induction process by obtaining the case $m$ from
the cases $m-2$ and $2$ we obtain even smaller constants. For example,
\begin{equation}
C_{\mathbb{R},m}=2^{\frac{1}{2}}\left(  \frac{C_{\mathbb{R},m-2}}
{A_{\frac{2m-4}{m-1}}^{2}}\right)  ^{\frac{m-2}{m}}\text{ for }m>3.
\end{equation}
In particular, if $2\leq m\leq14$ a careful calculation gives us
\begin{equation}
C_{\mathbb{R},m}=2^{\frac{m^{2}+6m-8}{8m}}\text{ if }m\text{ is even,
and}\label{unaestrella1}%
\end{equation}
\begin{equation}
C_{\mathbb{R},m}=2^{\frac{m^{2}+6m-7}{8m}}\text{ if }m\text{ is odd.}%
\label{unaestrella2}%
\end{equation}
However, in the next theorem a different induction approach leads us to even
smaller (and, thus, sharper) constants.

\begin{remark}
It is worth mentioning that the values for the above constants are not
explicitly calculated in \cite{defant} but, of course, are derived from
\cite{defant}. Since our procedure below will improve the constants, we will
not give much detail on the above estimates.
\end{remark}

The paper \cite{defant} provides, in fact, a family of constants $C_{m}$ for
Bohnenblust--Hille inequality. In this sense the following result is
essentially contained in \cite{defant}. However, for the sake of completeness
we prefer to sketch the proof from \cite{defant}. Our particular approach is
chosen to obtain sharper constants:

\begin{theorem}
\label{2_2} For every positive integer $m$ and real Banach spaces
$X_{1},\ldots,X_{m},$
\[
\Pi_{(\frac{2m}{m+1};1)}(X_{1},\ldots,X_{m};\mathbb{R})=\mathcal{L}%
(X_{1},...,X_{m};\mathbb{R})\text{ and }\left\Vert .\right\Vert
_{\pi(\frac{2m}{m+1};1)}\leq C_{\mathbb{R},m}\left\Vert .\right\Vert
\]
with
\[
C_{\mathbb{R},2}=2^{\frac{1}{2}}\text{ and }C_{\mathbb{R},3}=2^{\frac{5}{6}},
\]
\[
C_{\mathbb{R},m}=\frac{C_{\mathbb{R},m/2}}{A_{\frac{2m}{m+2}}^{m/2}}%
\]
for $m$ even and
\[
C_{\mathbb{R},m}=\left(  \frac{C_{\mathbb{R},\frac{m-1}{2}}}{A_{\frac
{2m-2}{m+1}}^{\frac{m+1}{2}}}\right)  ^{f(\frac{2m-2}{m+1},\frac{2m+2}{m+3}%
)}.\left(  \frac{C_{\mathbb{R},\frac{m+1}{2}}}{A_{\frac{2m+2}{m+3}}%
^{\frac{m-1}{2}}}\right)  ^{f(\frac{2m+2}{m+3},\frac{2m-2}{m+1})}%
\]
for $m$ odd. In particular,
\[
C_{\mathbb{R},m}=2^{\frac{1}{2}}C_{\mathbb{R},m/2}%
\]
if $m$ is even and $2\leq m\leq24.$
\end{theorem}

\begin{proof}
We start with the case $m$ even. The cases $m=2$ (Littlewood's $4/3$%
-inequality) and $m=3$ are known. We proceed by induction, obtaining the case
$m$ as a combination of the cases $m/2$ and $m/2$ (instead of $1$ and $m-1$ as
in \cite[Corollary 5.2]{defant}).\newline

\noindent Suppose that the result is true for $m/2$ and let us prove for $m$.
Let $U\in\mathcal{L}(X_{1},\ldots,X_{m};\mathbb{R})$ and $N$ be a positive
integer. For each $1\leq k\leq m$ consider $x_{1}^{(k)},\ldots,x_{N}^{(k)}\in
X_{k}$ so that $\left\Vert (x_{j}^{(k)})_{j=1}^{N}\right\Vert _{w,1}\leq1,$
$k=1,..,m.$

Consider, following the notation of Theorem \ref{b},
\[
q=2,\; s_{1}=s_{2}=\frac{2 \cdot(\frac{m}{2})}{\frac{m}{2}+1}=\frac{2m}{m+2}.
\]
Thus,
\[
w(s_{1},s_{2})=\frac{2m}{m+1}%
\]
and, from Theorem \ref{b}, we have
\begin{align*}
&  \displaystyle\left(  \sum\limits_{i_{1},\ldots,i_{m}=1}^{N}\left\vert
U(x_{i_{1}}^{(1)},\ldots,x_{i_{m}}^{(m)})\right\vert ^{\frac{2m}{m+1}}\right)
^{(m+1)/2m}\leq\\
&  \displaystyle\leq\left(  \sum\limits_{i_{1},\ldots,i_{m/2}=1}^{N}\left\Vert
\left(  U(x_{i_{1}}^{(1)},\ldots,x_{i_{m}}^{(m)})\right)  _{i_{\left(
m/2\right)  +1},\ldots,i_{m}=1}^{N}\right\Vert _{2}^{s_{2}}\right)
^{f(s_{2},s_{1})/s_{2}}\\
&  \left(  \sum\limits_{i_{\left(  m/2\right)  +1},\ldots,i_{m}=1}%
^{N}\left\Vert \left(  U(x_{i_{1}}^{(1)},\ldots,x_{i_{m}}^{(m)})\right)
_{i_{1}, \ldots,i_{m/2} =1}^{N}\right\Vert _{2}^{s_{1}}\right)  ^{f(s_{1}%
,s_{2})/s_{1}}.
\end{align*}

Now we need to estimate the two factors above. We will write $dt:=dt_{1}\ldots
dt_{m/2}$. For each $i_{\left(  m/2\right)  +1},...,i_{m}$ fixed, we have
(from Theorem \ref{d}),
\begin{align*}
&  \displaystyle\left\Vert \left(  U(x_{i_{1}}^{(1)},\ldots,x_{i_{m}}%
^{(m)})\right)  _{i_{1},\ldots,i_{m/2}=1}^{N}\right\Vert _{2}^{s_{1}}\leq\\
&  \displaystyle\leq\left(  A_{2,s_{1}}^{m/2}\right)  ^{s_{1}}{\textstyle}%
{\displaystyle\int\limits_{[0,1]^{m/2}}}\left\vert \sum\limits_{i_{1}%
,\ldots,i_{m/2}=1}^{N}r_{i_{1}}(t_{1})\ldots r_{i_{m/2}}(t_{m/2})U(x_{i_{1}%
}^{(1)},\ldots,x_{i_{m}}^{(m)})\right\vert ^{s_{1}}dt\\
&  \displaystyle=\left(  A_{2,s_{1}}^{m/2}\right)  ^{s_{1}}{\textstyle}%
{\displaystyle\int\limits_{[0,1]^{m/2}}}\displaystyle\left\vert U\left(
\sum\limits_{i_{1}=1}^{N}r_{i_{1}}(t_{1})x_{i_{1}}^{(1)},\ldots,\sum
\limits_{i_{m/2}=1}^{N}r_{i_{m/2}}(t_{m/2})x_{i_{m/2}}^{(m/2)},x_{i_{\left(
m/2\right)  +1}}^{(m/2)+1}\ldots,x_{i_{m-1}}^{(m-1)},x_{i_{m}}^{(m)}\right)
\right\vert ^{s_{1}}dt.
\end{align*}
Summing over all $i_{\left(  m/2\right)  +1},\ldots,i_{m}=1,\ldots,N$ we
obtain\newline$\displaystyle\sum\limits_{i_{\left(  m/2\right)  +1}%
,\ldots,i_{m}=1}^{N}\left\Vert \left(  U(x_{i_{1}}^{(1)},\ldots,x_{i_{m}%
}^{(m)})\right)  _{i_{1},\ldots,i_{m/2}=1}^{N}\right\Vert _{2}^{s_{1}}\leq
$\newline$\displaystyle\left(  A_{2,s_{1}}^{m/2}\right)  ^{s_{1}}%
{\textstyle}{\displaystyle\int\limits_{[0,1]^{m/2}}}\displaystyle\sum
\limits_{i_{\left(  m/2\right)  +1},\ldots,i_{m}=1}^{N}\left\vert U\left(
\sum\limits_{i_{1}=1}^{N}r_{i_{1}}(t_{1})x_{i_{1}}^{(1)},\ldots,\sum
\limits_{i_{m/2}=1}^{N}r_{i_{m/2}}(t_{m/2})x_{i_{m/2}}^{(m/2)},\ldots
,x_{i_{m}}^{(m-1)},x_{i_{m}}^{(m)}\right)  \right\vert ^{s_{1}}dt.$ Using the
case $m/2$ we thus have
\begin{align*}
&  \sum\limits_{i_{\left(  m/2\right)  +1},\ldots,i_{m}=1}^{N}\left\Vert
\left(  U(x_{i_{1}}^{(1)},\ldots,x_{i_{m}}^{(m)})\right)  _{i_{1}%
,\ldots,i_{m/2}=1}^{N}\right\Vert _{2}^{s_{1}}\leq\\
&  \leq\left(  A_{2,s_{1}}^{m/2}\right)  ^{s_{1}}{\displaystyle\int
\limits_{[0,1]^{m/2}}}\left\Vert U\left(  \sum\limits_{i_{1}=1}^{N}r_{i_{1}%
}(t_{1})x_{i_{1}}^{(1)},\ldots,\sum\limits_{i_{m/2}=1}^{N}r_{i_{m/2}}%
(t_{m/2})x_{i_{m/2}}^{(m/2)},\ldots\right)  \right\Vert _{\pi(s_{1};1)}%
^{s_{1}}dt\\
&  \leq\left(  A_{2,s_{1}}^{m/2}\right)  ^{s_{1}}\left(  \left\Vert
U\right\Vert C_{\mathbb{R},m/2}\right)  ^{s_{1}}.
\end{align*}
Hence
\[
\left(  \sum\limits_{i_{\left(  m/2\right)  +1},\ldots,i_{m}=1}^{N}\left\Vert
\left(  U(x_{i_{1}}^{(1)},\ldots,x_{i_{m}}^{(m)})\right)  _{i_{1}%
,\ldots,i_{m/2}=1}^{N}\right\Vert _{2}^{s_{1}}\right)  ^{\frac{1}{s_{1}}}\leq
A_{2,s_{1}}^{m/2}\left\Vert U\right\Vert C_{\mathbb{R},m/2}.
\]

The other estimate is exactly the same. Hence, combining both estimates, we
obtain
\[
\left(  \sum\limits_{i_{1},\ldots,i_{m}=1}^{N}\left\vert U(x_{i_{1}}%
^{(1)},\ldots,x_{i_{m}}^{(m)})\right\vert ^{\frac{2m}{m+1}}\right)
^{(m+1)/2m}\leq\left[  A_{2,s_{1}}^{m/2}\left\Vert U\right\Vert C_{\mathbb{R}%
,m/2}\right]  ^{1/2}\left[  A_{2,s_{1}}^{m/2}\left\Vert U\right\Vert
C_{\mathbb{R},m/2}\right] ^{1/2}.
\]
and
\begin{align*}
\left(  \sum\limits_{i_{1},\ldots,i_{m}=1}^{N}\left\vert U(x_{i_{1}}
^{(1)},\ldots,x_{i_{m}}^{(m)})\right\vert ^{\frac{2m}{m+1}}\right)
^{(m+1)/2m}  &  \leq A_{2,s_{1}}^{m/2}\left\Vert U\right\Vert C_{\mathbb{R}%
,m/2}\\
&  =A_{2,\frac{2m}{m+2}}^{m/2}\left\Vert U\right\Vert C_{\mathbb{R},m/2}.
\end{align*}
Hence
\[
C_{\mathbb{R},m}=A_{2,\frac{2m}{m+2}}^{m/2}C_{\mathbb{R},m/2}.
\]

Now let us estimate the constants $C_{\mathbb{R},m}.$ We know that $B_{2}=1$
and, from \cite{haag}, we also know that $A_{p}=2^{\frac{1}{2}-\frac{1}{p}}%
$whenever $p\leq1.847.$ So, for $2\leq m\leq24$ we have%
\begin{equation}
A_{\frac{2m}{m+2}}=2^{\frac{1}{2}-\frac{m+2}{2m}}=2^{\frac{-1}{m}%
}.\label{dosestrellas}%
\end{equation}

Hence, from (\ref{kcc}) and using the best constants of Khinchine's inequality
from \cite{haag}, we have%

\[
A_{2,\frac{2m}{m+2}}\leq A_{\frac{2m}{m+2}}^{-1}=2^{\frac{1}{m}}%
\]

and%
\[
C_{\mathbb{R},m}\leq\left(  2^{\frac{1}{m}}\right)  ^{m/2}C_{\mathbb{R}%
,m/2}=2^{\frac{1}{2}}C_{\mathbb{R},m/2}%
\]

for $m$ even, $2\leq m\leq24.$

The numerical values of $C_{\mathbb{R},m}$, for $m>24$, can be easily
calculated by using the exact values of $A_{\frac{2m}{m+2}}$ (see
\cite{haag}):
\[
A_{\frac{2m}{m+2}}=\sqrt{2}\left(  \frac{\Gamma\left(  \frac{\frac{2m}{m+2}%
+1}{2}\right)  }{\sqrt{\pi}}\right)  ^{(m+2)/2m}.
\]

For the case $m$ odd we proceed by induction, but the case $m$ is obtained as
a combination of the cases $\frac{m-1}{2}$ with $\frac{m+1}{2}$ instead of $1$
and $m-1$ as in \cite[Corollary 5.2]{defant}.

Consider, in the notation of Theorem \ref{b},
\[
q=2,\;s_{1}=\frac{2\left(  \frac{m-1}{2}\right)  }{\frac{m-1}{2}+1}%
=\frac{2m-2}{m+1}\text{ and }s_{2}=\frac{2\left(  \frac{m+1}{2}\right)
}{\frac{m+1}{2}+1}=\frac{2m+2}{m+3}.
\]
Thus,%
\[
w(s_{1},s_{2})=\frac{2m}{m+1}%
\]
and a similar proof gives us%
\begin{align*}
C_{\mathbb{R},m}  &  =\left(  A_{2,\frac{2m-2}{m+1}}^{\frac{m+1}{2}%
}C_{\mathbb{R},\frac{m-1}{2}}\right)  ^{f(\frac{2m-2}{m+1},\frac{2m+2}{m+3}%
)}.\left(  A_{2,\frac{2m+2}{m+3}}^{\frac{m-1}{2}}C_{\mathbb{R},\frac{m+1}{2}%
}\right)  ^{f(\frac{2m+2}{m+3},\frac{2m-2}{m+1})}\\
&  \leq\left(  \frac{C_{\mathbb{R},\frac{m-1}{2}}}{A_{\frac{2m-2}{m+1}}%
^{\frac{m+1}{2}}}\right)  ^{f(\frac{2m-2}{m+1},\frac{2m+2}{m+3})}.\left(
\frac{C_{\mathbb{R},\frac{m+1}{2}}}{A_{\frac{2m+2}{m+3}}^{\frac{m-1}{2}}%
}\right)  ^{f(\frac{2m+2}{m+3},\frac{2m-2}{m+1})}%
\end{align*}

\end{proof}

In the below table we compare the first constants $C_{m}=2^{\frac{m-1}{2}}$
and the constants that can be derived from \cite[Cor. 5.2]{defant} with the
\emph{new} constants $C_{\mathbb{R},m}$:

\begin{center}%
\begin{tabular}
[c]{c|c|c|c}%
$m$ & {\tiny New constants } - $C_{\mathbb{R},m}$ & $C_{m}${\tiny \ derived
from (\cite[Corollary 5.2]{defant}, 2010)} & $C_{m}=2^{\frac{m-1}{2}}$
\ {\tiny (\cite{Ka}, 1978)}\\\hline
$3$ & $\approx1.782$ & $2^{5/6}\approx1.782$ & $2^{2/2}=2$\\
$4$ & $2$ & $2^{\frac{18}{16}}\approx2.18$ & $2^{3/2}\approx2.828$\\
$5$ & $\approx2.298$ & $2^{\frac{28}{20}}\approx2.639$ & $2^{2}=4$\\
$6$ & $\approx2.520$ & $2^{\frac{40}{24}}\approx3.17$ & $2^{5/2}\approx
5.656$\\
$7$ & $\approx2.6918$ & $2^{\frac{54}{28}}\approx3.807$ & $2^{6/2}=8$\\
$8$ & $\approx2.8284$ & $2^{\frac{70}{32}}\approx4.555$ & $2^{7/2}%
\approx11.313$\\
$9$ & $\approx3.055$ & $2^{\frac{88}{36}}\approx5.443$ & $2^{8/2}=16$\\
$10$ & $\approx3.249$ & $2^{\frac{108}{40}}\approx6.498$ & $2^{9/2}%
\approx22.627$\\
$11$ & $\approx3.4174$ & $2^{\frac{130}{44}}\approx7.752$ & $2^{10/2}=32$\\
$12$ & $\approx3.563$ & $2^{\frac{154}{48}}\approx9.243$ & $2^{11/2}%
\approx45.254$%
\end{tabular}

\end{center}

\section{Improved constants for the Bohnenblust--Hille theorem: The complex
case}

As in the real case, following the proof of \cite[Corollary 5.2]{defant} and
using the optimal values for the constants of Khinchine's inequality (due to
Haagerup) and using that $K_{G}=C_{\mathbb{C},2}$ (see \cite{monats}), the
following estimates can be calculated for $C_{m}$:
\begin{align*}
C_{\mathbb{C},2}  &  = K_{G}\leq1.4049<\sqrt{2},\\
C_{\mathbb{C},m}  &  = 2^{\frac{m-1}{2m}}\left(  \frac{C_{\mathbb{C},m-1}%
}{A_{\frac{2m-2}{m}}}\right)  ^{1-\frac{1}{m}}\text{ for }m\geq3,
\end{align*}
In particular, if $2\leq m\leq13$,
\begin{equation}
C_{\mathbb{C},m}\leq2^{\frac{m^{2}+m-6}{4m}} K_{G}^{2/m} .\nonumber
\end{equation}

The above estimates improve the values $C_{\mathbb{C},m}=2^{\frac{m-1}{2}}$
but are worst than the constants $C_{\mathbb{C},m}=\left( \frac{2}{\sqrt{\pi}%
}\right)  ^{m-1}$ obtained by Defant and Sevilla-Peris \cite{defant2}.
However, our approach will provide even better constants.

The following lemma is essentially the main result from the previous section
which comes from \cite{defant}, although now we will obtain different
constants, since we will be dealing with the complex case.

\begin{lemma}
For every positive integer $m$ and complex Banach spaces $X_{1},\dots,X_{m}$,
\[
\Pi_{(\frac{2m}{m+1};1)}(X_{1},\ldots,X_{m};\mathbb{C})=\mathcal{L}
(X_{1},\ldots,X_{m};\mathbb{C})\text{ and }\left\Vert .\right\Vert _{\pi
(\frac{2m}{m+1};1)}\leq C_{\mathbb{C},m}\left\Vert .\right\Vert
\]
with
\[
C_{\mathbb{C},m}=\left(  \frac{2}{\sqrt{\pi}}\right)  ^{m-1}\text{ for }m=2,3,
\]

\[
C_{\mathbb{C},m}=\frac{C_{\mathbb{C},m/2}}{A_{\frac{2m}{m+2}}^{m/2}}%
\]
for $m$ even and
\[
C_{\mathbb{C},m}=\left(  \frac{C_{\mathbb{C},\frac{m-1}{2}}}{A_{\frac
{2m-2}{m+1}}^{\frac{m+1}{2}}}\right)  ^{f(\frac{2m-2}{m+1},\frac{2m+2}{m+3}%
)}\left(  \frac{C_{\mathbb{C},\frac{m+1}{2}}}{A_{\frac{2m+2}{m+3}}^{\frac
{m-1}{2}}}\right) ^{f(\frac{2m+2}{m+3},\frac{2m-2}{m+1})}%
\]
for $m$ odd.
\end{lemma}

The proof of this result is essentially in the same spirit as that of Theorem
\ref{2_2}. The cases of $C_{\mathbb{C},2}$ and $C_{\mathbb{C},3}$ are already
known and the proof is (also) done by induction.

\subsection{Comparing the ``\emph{first}'' constants}

The first constants $D_{m}=\left(  \frac{2}{\sqrt{\pi}}\right) ^{m-1}$ from
\cite{defant2} are better than the constants that we have obtained in the
previous lemma. However our constants present a smaller asymptotical growth
and, from a certain level $m$ on, they are better than $D_{m}=\left(  \frac
{2}{\sqrt{\pi}}\right) ^{m-1}.$ As we see next, this occurs when $m\geq
7$.\newline

\noindent Here below we compare the first constants. For the case $m=4$,
notice that we have
\[
C_{\mathbb{C},4}=2^{1/2} \cdot C_{\mathbb{C},2}=2^{1/2} \cdot\left(  \frac
{2}{\sqrt{\pi}}\right)  \approx1.5957
\]
but this constant is \emph{worst} than $\left(  \frac{2}{\sqrt{\pi}}\right)
^{3}\approx1.437.$ So, in order to improve the constants that follow, it would
be better to consider $\left(  \frac{2}{\sqrt{\pi}}\right) ^{3}$ instead of
$2^{1/2}\cdot C_{\mathbb{C},2}$ for the value of $C_{\mathbb{C},4}$.

Similarly, for the cases $5\leq m\leq6$ we also have that $C_{\mathbb{C},m}$
is slightly \emph{worst} than $(2/\sqrt{\pi})^{m-1}$ but, for $m=7$ we have
\[
C_{\mathbb{C},7}=\left(  \frac{\left(  \frac{2}{\sqrt{\pi}}\right) ^{2}%
}{\left(  2^{\frac{1}{2}-\frac{8}{12}}\right) ^{4}}\right) ^{3/7} \cdot\left(
\frac{\left(  \frac{2}{\sqrt{\pi}}\right) ^{3}}{\left( 2^{\frac{1}{2}%
-\frac{10}{16}}\right) ^{3}}\right) ^{4/7}=1.9293<\left(  \frac{2}{\sqrt{\pi}%
}\right) ^{6}%
\]
and our constants are \emph{better} than the \emph{old} ones. Also, as we
announced, for $m\geq 7$ our constants also improve the \emph{old} ones.\newline

\noindent In the next section we state the previous lemma using the
information we just obtained.

\subsection{Comparing the \emph{\textquotedblleft remaining\textquotedblright}
constants ($m>7$).}

Now it is time to state the last lemma adding the \emph{better} constants:

\begin{theorem}
For every positive integer $m$ and every complex Banach spaces $X_{1},\ldots,
X_{m},$
\[
\Pi_{(\frac{2m}{m+1};1)}(X_{1}, \ldots,X_{m};\mathbb{C})=\mathcal{L}(X_{1},
\ldots,X_{m};\mathbb{C})\text{ and }\left\Vert .\right\Vert _{\pi(\frac
{2m}{m+1};1)}\leq C_{\mathbb{C},m}\left\Vert .\right\Vert
\]
with

\begin{itemize}
\item $\displaystyle C_{\mathbb{C},m}=\left(  \frac{2}{\sqrt{\pi}}\right)
^{m-1}\text{ for } m \in\{2,3,4,5,6\},$

\item $\displaystyle C_{\mathbb{C},m}=\frac{C_{\mathbb{C},m/2}}{A_{\frac
{2m}{m+2}}^{m/2}}\text{ }$ for $m > 6$ even, and

\item $\displaystyle C_{\mathbb{C},m}=\left(  \frac{C_{\mathbb{C},\frac
{m-1}{2}}}{A_{\frac{2m-2}{m+1}}^{\frac{m+1}{2}}}\right) ^{f(\frac{2m-2}%
{m+1},\frac{2m+2}{m+3})} \cdot\left( \frac{C_{\mathbb{C},\frac{m+1}{2}}%
}{A_{\frac{2m+2}{m+3}}^{\frac{m-1}{2}}}\right) ^{f(\frac{2m+2}{m+3}%
,\frac{2m-2}{m+1})}$ for $m>5$ odd.
\end{itemize}
\end{theorem}

The following table compares these new constants with the previous ones:

\begin{center}%
\begin{tabular}
[c]{c|c|c|c|c}%
$m$ & {\tiny New Constants} - $C_{\mathbb{C},m}$ & $\left(  \frac{2}{\sqrt
{\pi}}\right)  ^{m-1}${\tiny (\cite{Que, defant2}, 1995)} & $2^{\frac{m-1}{2}%
}${\tiny (\cite{Ka},1978)} & $m^{\frac{m+1}{2m}}2^{\frac{m-1}{2}}%
${\tiny (\cite{bh},1931)}\\\hline
$8$ & $\approx2.031$ & $\approx2.329$ & $\approx11.313$ & $\approx36.442$\\
$9$ & $\approx2.172$ & $\approx2.628$ & $16$ & $\approx54.232$\\
$10$ & $\approx2.292$ & $\approx2.965$ & $\approx22.627$ & $\approx80.283$\\
$11$ & $\approx2.449$ & $\approx3.346$ & $32$ & $\approx118.354$\\
$12$ & $\approx2.587$ & $\approx3.775$ & $\approx45.425$ & $\approx173.869$\\
$13$ & $\approx2.662$ & $\approx4.260$ & $64$ & $\approx254.680$\\
$14$ & $\approx2.728$ & $\approx4.807$ & $\approx90.509$ & $\approx372.128$\\
$15$ & $\approx2.805$ & $\approx5.425$ & $128$ & $\approx542.574$\\
$16$ & $\approx2.873$ & $\approx6.121$ & $\approx181.019$ & $\approx789.612$%
\end{tabular}

\end{center}

\bigskip

\noindent\textbf{Acknowledgements.} The authors thank A. Defant and P.
Sevilla-Peris for important remarks on the complex case. The authors are also
very grateful to the referee, whose thorough analysis and insightful remarks
helped to obtain even sharper constants.

\end{document}